\documentclass[a4paper]{article}
\usepackage{amsmath}\usepackage{epsf,amsfonts,amsthm}
\usepackage{fontenc,indentfirst, delarray,amsfonts,amsmath,amssymb}

\newcommand{\be}{\begin{equation}}
\newcommand{\ee}{\end{equation}}
\newcommand{\bea}{\begin{eqnarray*}}
\newcommand{\eea}{\end{eqnarray*}}

\newcommand{\Ci}{C^{\infty}}

\newcommand{\p}{\partial}

\newcommand{\cS}{{\cal S}}

\newcommand{\raa}{\rightarrow}

\newcommand{\E}{\ell}
\newcommand{\m}{\!\!\mid}

\newcommand{\N}{{\mathbb{N}}}
\newcommand{\Z}{{\mathbb{Z}}}
\newcommand{\R}{{\mathbb{R}}}

\newcommand{\lp}{\left(}
\newcommand{\rp}{\right)}

\newcommand{\ap}{\alpha}

\newcommand{\op}[1]{\!\!\mathop{\rm ~#1}\nolimits}
\newcommand{\mb}{\mathbb}

\newcommand{\cD}{{\cal D}}

\newcommand{\cX}{{\cal X}}

\mathchardef\za="710B  %\alpha
\mathchardef\zb="710C  %\beta
\mathchardef\zg="710D  %\gamma
\mathchardef\zd="710E  %\delta
\mathchardef\zve="710F %\epsilon
\mathchardef\zz="7110  %\zeta
\mathchardef\zh="7111  %\eta
\mathchardef\zvy="7112 %\theta
\mathchardef\zi="7113  %\iota
\mathchardef\zk="7114  %\kappa
\mathchardef\zl="7115  %\lambda
\mathchardef\zm="7116  %\mu
\mathchardef\zn="7117  %\nu
\mathchardef\zx="7118  %\xi
\mathchardef\zp="7119  %\pi
\mathchardef\zr="711A  %\rho
\mathchardef\zs="711B  %\sigma
\mathchardef\zt="711C  %\tau
\mathchardef\zu="711D  %\upsilon
\mathchardef\zvf="711E %\phi
\mathchardef\zq="711F  %\chi
\mathchardef\zc="7120  %\psi
\mathchardef\zw="7121  %\omega
\mathchardef\ze="7122  %\varepsilon
\mathchardef\zy="7123  %\vartheta
\mathchardef\zf="7124  %\varomega
\mathchardef\zvr="7125 %\varrho
\mathchardef\zvs="7126 %\varsigma
\mathchardef\zf="7127  %\varphi
\mathchardef\zG="7000  %\Gamma
\mathchardef\zD="7001  %\Delta
\mathchardef\zY="7002  %\Theta
\mathchardef\zL="7003  %\Lambda
\mathchardef\zX="7004  %\Xi
\mathchardef\zP="7005  %\Pi
\mathchardef\zS="7006  %\Sigma
\mathchardef\zU="7007  %\Upsilon
\mathchardef\zF="7008  %\Phi
\mathchardef\zW="700A  %\Omega

\newcommand{\beas}{\begin{eqnarray*}}
\newcommand{\eeas}{\end{eqnarray*}}

\begin{document}
\title{On quantum and classical\\ Poisson algebras}\author{Janusz GRABOWSKI\thanks{Research
of Janusz Grabowski supported by the Polish Ministry of Scientific
Research and Information Technology under the grant No. 2 P03A
02024.}\ , Norbert PONCIN\thanks{The research of Norbert Poncin
was supported by grant R1F105L10.}} \maketitle

\newtheorem{rem}{Remark}
\newtheorem{theo}{Theorem}
\newtheorem{lem}{Lemma}
\newtheorem{prop}{Proposition}
\newtheorem{cor}{Corollary}

\begin{abstract} Results on derivations and automorphisms of some
quantum and classical Poisson algebras,
as well as characterizations of manifolds by the Lie structure of
such algebras, are revisited and extended. We prove in particular
somehow unexpected fact that the algebras of linear differential
operators acting on smooth sections of two real vector bundles of
rank 1 are isomorphic as Lie algebras if and only if the base
manifolds are diffeomorphic, independently whether the line
bundles themselves are isomorphic or not.
\end{abstract}

\bigskip\noindent
\textit{\textbf{MSC 2000:} 17B63 (Primary), 13N10, 16S32, 17B40,
17B65, 53D17 (Secondary).}

\medskip\noindent
\textit{\textbf{Key words:} line bundles, differential operators,
derivations, automorphisms, Lie algebras, Chevalley cohomology,
one-parameter groups.}

\section{Introduction}

Let us start with an overview of relevant literature on
isomorphisms and derivations of infinite-dimensional Lie algebras.
\\

In \cite{PS}, Pursell and Shanks proved the well-known result
stating that the Lie algebra of all smooth compactly supported
vector fields of a smooth manifold characterizes the
differentiable structure of the variety. Similar upshots were
obtained in numerous subsequent papers dealing with different Lie
algebras of vector fields and related algebras (see e.g.
\cite{A,Am,AG,JG,JG3,HM,O,S}).

Derivations of certain infinite-dimensional Lie algebras arising
in Geometry were also studied in different situations (note that
in infinite dimension there is no such a clear correspondence
between derivations and one-parameter groups of automorphisms as
in the finite-dimensional case). Let us mention a result of
L.~S.~Wollenberg \cite{Wol} who described all derivations of the
Lie algebra of polynomial functions on the canonical symplectic
space $\R^2$ with respect to the Poisson bracket. It turned out
that there are outer derivations of this algebra in contrast to
the corresponding Weyl algebra. This can be viewed as a variant of
a "no-go" theorem (see \cite{J}) stating that the Dirac
quantization problem \cite{Dir} cannot be solved satisfactorily
because the classical and the corresponding quantum algebras are
not isomorphic as Lie algebras. An algebraic generalization of the
latter fact, known as the {\it algebraic "no-go" theorem,} has
been proved in \cite{GG} by different methods. Derivations of the
Poisson bracket of all smooth functions on a symplectic manifold
have been determined in \cite{ADML} (for the real-analytic case,
see \cite{Gr2}). Another important result is the one by F.~Takens
\cite{Tak} stating that all derivations of the Lie algebra
$\cX(M)$ of smooth vector fields on a manifold $M$ are inner. The
same turned out to be valid for analytic cases \cite{Gr1}. Some
cases of the Lie algebras of vector fields associated with
different geometric structures were studied in a series of papers
by Y.~Kanie \cite{Ka1}--\cite{Ka4}.

Our work \cite{GP} contains Shanks-Pursell type results for the
Lie algebra $\cD(M)$ of all linear differential operators of a
smooth manifold $M$, for its Lie subalgebra $\cD^1(M)$ of all
linear first-order differential operators of $M$, and for the
Poisson algebra $\cS(M)=\op{Pol}(T^*M)$ of all polynomial
functions on $T^*M,$ the symbols of the operators in $\cD(M).$
Furthermore, we computed all the automorphisms of these algebras
and showed that the Lie algebras $\cD(M)$ and $\cS(M)$ are not
integrable. The paper \cite{GP3} provides their derivations, so it
is a natural continuation of this previous work and can be
considered as a generalization of the results of Wollenberg and
Takens. It is also shown which derivations generate one-parameter
groups of automorphisms and the explicit form of such
one-parameter groups
is given.\\

The first part of the present text is an intuitive description of
the major facts explained in \cite{GP3}. Moreover, experience of
different approaches to the quantization problem and the geometric
study of differential equations incite to substitute differential
operators acting on tensor densities for differential operators on
functions. In the frame of our previous works, this substitution
requires investigations on a possible characterization of a
manifold $M$ by the Lie algebra of differential operators acting
on densities on $M$ of arbitrary fixed weight, or more generally,
on the potential characterization, by the canonical Lie algebra
structure of the space of linear differential operators on smooth
sections of an arbitrary $\mathbb{R}$-line bundle $L$, of the base
manifold $M$ or even of the bundle $L$ itself. These problems are
solved in the second part of this paper.

\section{Derivations of some quantum and classical Poisson algebras}

\subsection{Locality and weight}

In this section we depict the derivations of the algebras ${\cal
D}^1(M)$, ${\cal S}(M)$, and ${\cal D}(M)$. Let $({\cal D},[.,.])$
be one of these three filtered Lie algebras and let $C$ be a
derivation of $({\cal D},[.,.])$. We speak of operators when
referring to elements of ${\cal D}$ and denote the algebra
$\Ci(M)$ of smooth functions of $M$ by ${\cal A}$. The adjoint
action of a smooth function $f$ of $M$, regarded as a differential
operator of order 0, on an operator $D\in{\cal D}^{i}$ lowers the
filtration degree by $1$. This provides a tool for proofs by
induction. The idea is really fruitful if derivations have weight
$0$. Indeed, a bracket such as $[CD,f]$, which involves the chosen
derivation, is then also a member of ${\cal D}^{i-1}$.

We first prove that any derivation $C$ has a bounded weight, i.e.
that there is a positive integer $d$, such that $C{\cal
D}^{i}\subset {\cal D}^{i+d},\forall i\in\N$. The proof uses the
derivation property on functions, the characterization of filters
``\`a la Vinogradov'', and the result that the ${\cal A}$-module
$\zW^1(M)$ of differential $1$-forms is spanned by the
differentials of a finite number of functions. This last upshot is
a consequence of the Whitney's embedding theorem.

In order to verify that investigation by local computations is
possible, we have to check if any derivation can be restricted to
a domain of local coordinates. This means that we must prove that
a derivation is always a local operator. We obtain locality using
a general technique worked out by De Wilde and Lecomte, see
\cite{DWLc1}: if an operator $D\in{\cal D}$ vanishes in a
neighborhood $U$ of a point $x\in M$, it reads
$D=\sum_k[X_k,D_k]$, where the sum is finite and the vector fields
$X_k$ and operators $D_k$ vanish in some neighborhood $V\subset U$
of $x$. The derivation property $CD=\sum_k\lp
[CX_k,D_k]+[X_k,CD_k]\rp$ then allows to conclude.

Let us emphasize that significant information on automorphisms and
derivations is encoded in the automorphism and derivation
properties written for two functions. If $(x^1,\ldots,x^n)$ are
local coordinates in an open subset $U\subset M$, we get
\[0=C[x^{i},x^j]=[Cx^{i},x^j]+[x^{i},Cx^j].\] The values $Cx^{i}$
are differential operators over $U$ or polynomials of $T^*U$. In
the first case, we symbolically write the derivatives in these
operators $Cx^{i}$ as monomials in the corresponding components
$(\zx_1,\dots,\zx_n)$ of some linear form
$\zx\in(\mathbb{R}^n)^*$. Hence $Cx^{i}\simeq P^{i}$, where the
$P^{i}$ are polynomials of $T^*U.$ In this polynomial language,
the above result reads \[\p_{\zx_j}P^{i}=\p_{\zx_i}P^{j}.\]
Integration furnishes us with a polynomial $P$ such that
$\p_{\zx_i} P=P^{i}$, i.e. $[P,x^{i}]=Cx^{i}$. So derivation $C$
coincides on coordinate functions with an interior derivation. It
is easily seen that for an arbitrary function $f$, the derivations
$C$ and $\op{ad}P$ differ by a function, $Cf-[P,f]\in{\cal A}$,
i.e. locally $C-\op{ad}P$ respects the lowest filter. After a
gluing process and a generalization to higher order filters, we
conclude that any derivation can be corrected by an interior
derivation in such a way that the filtration is respected. We will
refer to this property as ``property P1''. Let us stress that
$P\in{\cal D}^{d+1}$ since the weight of $C$ is $d$. Operator $P$
is not unique, the set of all convenient $P$ is $P+{\cal D}^1$.

\subsection{Restriction to functions\label{corrections}}

We prove that for any derivation $C$ that respects the filtration
there is a unique vector field $Y$, such that the derivation
$C-\op{ad}Y$ respects the filtration and reduces on functions to a
multiple of identity, i.e. $\lp C-\op{ad}Y\rp\m_{\cal
A}=\kappa\,\op{id},$ where $\kappa\in\mathbb{R}$ is uniquely
determined by $C$.

In the following we refer to this result as ``property P2''.
Unless differently stated, we assume that all derivations examined
below have been corrected and have acquired properties P1 and P2.

The proof of the preceding upshot is based upon a technique
similar to that used in \cite[Sect. 2.3.2]{GP2} and will not be
described here.

\subsection{Derivations of first order linear differential operators\label{derfristdiffop}}

The Lie algebra of first order differential operators has a
canonical splitting, ${\cal D}^1(M)={\cal A}\oplus{\cal X}(M)$,
where ${\cal X}(M)$ is the Lie algebra of vector fields of $M$. In
the following we simply write ${\cal D}^1$ and ${\cal X}$, if no
misunderstanding is possible. In view of property P2, we have
$Cf=\kappa f$ ($f\in{\cal A},\kappa\in\mathbb{R}$). The derivation
property shows that $C\m_{\cal X}$ is a $1$-cocycle of the
canonical representation of the Lie algebra ${\cal X}$ on the
space ${\cal A}.$ These cocycles are well-known (see \cite{Fuc},
\cite{DWLc2}): $CX=\zl\op{div}X+\zw(X)$ ($X\in{\cal
X},\zl\in\mathbb{R},\zw\in\zW^1(M)\cap\op{ker}\op{d}$). So we know
$C$ on any first-order operator. If we wish to recover the
initial, not yet corrected (see \ref{corrections}) derivation (we
denote it also by $C$), we have to add again the corrections.
Finally,\[C(f+X)=[Y,f+X]+\kappa f+\zl\op{div}X+\zw(X),\forall
f\in{\cal A},\forall X\in{\cal X},\] where $Y\in{\cal
X},\kappa,\zl\in\mathbb{R},\zw\in\zW^1(M)\cap\op{ker\op{d}}$ are
uniquely defined by $C$. The cohomological translation of this
result is \[H^1({\cal D}^1,{\cal D}^1)=\mathbb{R}^2\oplus
H^1_{DR}(M).\] Here $H^1({\cal D}^1,{\cal D}^1)$ is the first
cohomology space of the Lie algebra ${\cal D}^1$ and $H^1_{DR}(M)$
is the first de Rham cohomology group of the underlying manifold
$M$.

Explanations regarding the divergence can be found in \cite[Sect.
2.5.1]{GP2}. Let us recall that any nowhere vanishing $1$-density
$\zr_0$ defines a vector space isomorphism $\zt_0$ between the
space of $1$-densities and the space of functions. Nevertheless
these spaces are not isomorphic as modules over the Lie algebra of
vector fields. Indeed, if ${\cal L}_X$ and $L_X$ denote the Lie
derivatives with respect to a vector field $X$, of $1$-densities
and functions respectively, the difference $\zt_0\circ{\cal
L}_X\circ\zt_0^{-1}-L_X$ is the value at $X$ of a $1$-cocycle of
${\cal X}$ with coefficients in ${\cal A}$. This cocycle is the
divergence implemented by $\zr_0$. There is no canonical
divergence. Nevertheless all divergences induced by nowhere
vanishing $1$-densities are cohomologous. So these divergences
define a privileged cohomology class. The divergence above and
below is a fixed divergence of this class.

\subsection{Derivations of polynomials on the cotangent bundle\label{derpoly}}

Any (corrected) derivation $C$ of ${\cal S}(M)=\op{Pol}(T^*M)$ (in
the following we simply note ${\cal S}$) restricts to a derivation
(still denoted by $C$) of the Lie algebra ${\cal S}^1$ of
polynomials of degree $1$ at most. Since this algebra is
isomorphic to ${\cal D}^1$, derivation $C$ reads
\begin{equation}C(f+X)=\kappa f+\zl\op{div}X+\zw(X),\forall f\in{\cal
A},\forall X\in{\cal X}.\label{derpolyrest}\end{equation} If we
impose the derivation condition, not only for elements of ${\cal
S}^1$ but for all polynomials in ${\cal S}$, the terms of the
r.h.s. of Equation (\ref{derpolyrest}) either cancel or turn out
to be the traces on the ${\cal S}^1$-level of derivations of the
whole algebra ${\cal S}$. In this intuitive approach we confine
ourselves to trying to extend these terms as derivations of ${\cal
S}$.

Since ${\cal S}$ is a graded algebra, one of its derivations is
the so-called {\it degree derivation}, \[\op{Deg} :{\cal S}_{i}\ni
P\raa (i-1)P\in{\cal S}_{i},\] which just multiplies by the
(shifted) degree of the argument. Of course ${\cal S}_i$ is the
space $\op{Pol}^{i}(T^*M)$ of homogeneous polynomials of degree
$i$. Visibly $-\kappa\op{Deg}$ is a derivation of ${\cal S}$ that
extends the first term of the r.h.s. of Equation
(\ref{derpolyrest}).

It should be clear that such an extension does not exist for the
second term $\zl\op{div}$.

Since $\zw$ is locally exact, its value at $X$ locally reads
\[\zw(X)=(\op{d}f)(X)=\{X,f\},\] where $f$ is a local function
and where $\{.,.\}$ is the standard Poisson bracket of $T^*M.$ So,
if, for any polynomial $P\in{\cal S}$, we locally define
\[\overline{\zw}(P):=\{P,f\}=\zL(\op{d}P,\op{d}f),\] where $\zL$
is the corresponding Poisson tensor, we see that $\overline{\zw}$
is a well and globally defined derivation of ${\cal S}$ that
extends our third term. It is obvious from the preceding equation
that $\overline{\zw}$ is the (vertical) vector field $\zw^{v}$ of
$T^*M$ induced by (the pullback of) $\zw$ (to the cotangent
bundle).

Finally we understand that any derivation $C$ of ${\cal S}$ is of
the type
\[C(P)=\{Q,P\}-\kappa\op{Deg}P+\zw^v(P),\forall P\in{\cal S},\]
where $Q\in{\cal S}, \zk\in\R, \zw\in\zW^1(M)\cap\op{ker}\op{d}$.
Let us still mention that $\zk$ is unique, whereas the set of
appropriate $(Q,\zw)$ is $\{(Q+h,\zw+\op{d}h), h\in{\cal A}\}$.
The cohomological version of this second result is \[H^1({\cal
S},{\cal S})=\R\oplus H_{\op{DR}}^1(M),\] with self-explaining
notations.

\subsection{Derivations of linear differential operators\label{derdiffop}}

The intuitive approach is as in Section \ref{derpoly}, conclusions
are similar. Note nevertheless that the degree derivation, which
extends the first term of Equation (\ref{derpolyrest}), is tightly
connected with the grading of the classical Poisson algebra ${\cal
S}$. Since the quantum algebra ${\cal D}(M)$ (${\cal D}$ for
short) is only filtered, we guess that such an extension is no
longer possible on the quantum level.

We now understand that any derivation $C$ of ${\cal D}$ has the
form
\[C(D)=[\zD,D]+\overline{\zw}(D),\forall D\in{\cal D},\]
where $\zD\in{\cal D},\zw\in\zW^1(M)\cap \op{ker}\op{d}.$ The
lowering (its weight with respect to the filtration degree is
$-1$) derivation $\overline{\zw}$ is defined as in Section
\ref{derpoly}. The convenient $(\zD,\zw)$ are again
$\{(\zD+h,\zw+\op{d}h),h\in{\cal A}\}$. Moreover,
\[H^1({\cal D},{\cal D})=H^1_{\op{DR}}(M).\]

\section{Canonical and equivariant quantizations}

Needless to say that the above depicted intuitive presentation is
merely a kind of ``story''. Proofs are completely different and
comparatively technical. So the quest for equations that
annihilate both constants $\zk$ and $\zl$ (see Eq.
(\ref{derpolyrest})) is, in the case of differential operators, a
little bit as the hunt for the proverbial needle in a haystack.
Computations can be reduced to the search for some operators that
intertwine the canonical action of vector fields. These
intertwining operators had been discovered in \cite{NP} and in
\cite{BHMP}.\\

The first of these papers as well as our work \cite{GP} are based
upon formal calculus. In the main, this symbolism consists in the
substitution of monomials $\xi_1^{\ap^1}\ldots\xi_n^{\ap^n}$ in
the components of a linear form $\xi\in (\R^n)^*$ to the
derivatives $\p_{x^1}^{\ap^1}\ldots\p_{x^n}^{\ap^n}f$ of a
function $f$. In other words, we exploit the so-called affine
symbol map,
\[\begin{array}{rl}\zs_{\op{aff}}:{\cal D}(\mathbb{R}^n)\ni D=\sum D^{i_1\ldots
i_j}(x)\!\!\!\!&\p_{x^{i_1}}\ldots\p_{x^{i_j}}\\
\rightarrow\zs_{\op{aff}}(D)=&\!\!\!\!\!\sum D^{i_1\ldots
i_j}(x)\zx_{i_1}\ldots \zx_{i_j}\in{\cal S}(\R^n)=
\op{Pol}(T^*\mathbb{R}^n),\end{array}\] where the differential
operator $D$ is decomposed in some coordinate system and where the
coefficients $D^{i_1\ldots i_j}(x)$ are smooth functions. Note
that in order to increase readability, we have explicitly written
the variables of the involved functions, thus identifying a
function with its value at a generic point. This affine symbol
method is known in Mechanics as the ``normal ordering'' or
``canonical symbolization/quantization''. Its systematic use in
Differential Geometry is originated in papers by M. Flato and A.
Lichnerowicz \cite{FL} as well as M. De Wilde and P. Lecomte
\cite{DWLc2}, dealing with the Chevalley-Eilenberg cohomology of
the Lie algebra of vector fields associated with the Lie
derivative of differential forms. This polynomial modus operandi
matured during the last twenty years and developed into a quite
powerful and universal computing technique that has been
successfully applied in numerous works (see e.g. \cite{LMT} or
\cite{NP2}). A more detailed explanation of the key features of
this non standard procedure can be found in \cite{NP}.\\

The reader might ask if tensor field $\zs_{\op{aff}}(D)$ is
independent on the chosen coordinates, i.e. if the vector space
isomorphism $\zs_{\op{aff}}$ commutes with diffeomorphisms. An
easy computation allows to convince oneself that only the highest
order terms, the principal symbol of $D$, have an intrinsic
meaning.

This problem is tightly connected with equivariant quantization.
Let us first replace, as already mentioned, differential operators
acting on functions by differential operators acting on densities
of arbitrary weight $\zl.$ Remember that $0$-densities are just
functions and that the space of $1/2$-densities is pre-Hilbert and
deserves particular attention. It seems natural to wonder whether
there is an isomorphism $\zs_{\chi}\neq\zs_{\op{aff}}$ between the
space ${\cal D}_{\zl}(\R^n)$ of differential operators on
densities of weight $\zl$ and the space $\op{Pol}(T^*\R^n)$ of
polynomials on the cotangent bundle, which commutes with
diffeomorphisms or---on the infinitesimal level---with the Lie
derivatives with respect to vector fields. The existence of such a
total symbol map would imply that the spaces of differential
operators on $\zl$- and $\mu$-densities are isomorphic as modules
over the Lie algebra of vector fields. The classification of these
modules has been obtained in \cite{LMT}. It shows that there is no
${\cal X}(\R^n)$-equivariant symbol map, that we have to relax the
invariance condition.

%Consider a general mechanical system with symmetries. When
%dividing these symmetries out, we pass to the reduced phase space
%that has some special structure, e.g. a projective structure. A
%manifold endowed with a projective structure is locally modelled
%on the projective space $\R P^n$ and transition functions are
%implemented by members of the special linear group
%$\op{SL}(n+1,\R)$. It is natural to require the equivariance
%condition only for these coordinate changes, i.e. for the action
%of $\op{SL}(n+1,\R)$.
Note however that the canonical action of $\op{SL}(n+1,\R)$ on
$\R^{n+1}$ induces
%not only an action on $\R P^n$, but also
a ``local'' action on $\R^n$, which in turn induces an action of
the Lie algebra $\op{sl}(n+1,\R)$ on $\R^n$. This tangent action
realizes $\op{sl}(n+1,\R)$ as a maximal Lie subalgebra
$\op{sl}_{n+1}$ of the algebra ${\cal X}_*(\R^n)$ of polynomial
vector fields on $\R^n.$ So the aforementioned relaxed
equivariance condition can be possibly changed into the
equivariance with respect to the action of this maximal algebra
$\op{sl}_{n+1}$ of infinitesimal projective transformations. More
precisely, we are looking for a vector space isomorphism
\[\zs_{\op{sl}}:{\cal
D}_{\zl}(\mathbb{R}^n)\rightarrow\op{Pol}(T^*\mathbb{R}^n),\] such
that \[\zs_{\op{sl}}\circ{\cal
L}_X=L_X\circ\zs_{\op{sl}},\;\forall X\in \op{sl}_{n+1}\] and
\[\zs_{\op{sl}}(D)-\zs (D)\in\op{Pol}^{\le
k-1}(T^*\mathbb{R}^n),\;\forall D\in{\cal
D}^k_{\zl}(\mathbb{R}^n),\] where ${\cal L}_X$ and $L_X$ are the
standard actions of a vector field $X$ on differential operators
and tensor fields respectively, where $\zs$ denotes the principal
symbol and where superscripts $\le k-1$ and $k$ are the filtration
degrees. The last condition is a normalization condition that
assures uniqueness of the projectively equivariant symbol and
quantization maps $\zs_{\op{sl}}$ and $\zs_{\op{sl}}^{-1}$.
Existence has also been proven, see \cite{PLVO}. Note that
$\zs_{\op{aff}}$ is called the affine symbol map since it
intertwines the Lie derivatives with respect to the affine vector
fields.

If $M$ is endowed with a projective structure, the projectively
equivariant quantization map $Q=\zs_{\op{sl}}^{-1}$, which exists
in any chart of any projective atlas, is of course well-defined on
$M$. This quantization procedure $Q:\op{Pol}(T^*M)\raa{\cal
D}_{\zl}(M)$ defines a family of invariant $*$-products on $T^*M$.
Indeed, it is easily checked that if for any formal parameter
$\hbar$ and any $P\in\op{Pol}^k(T^*M)$, we set
\[Q_{\hbar}P=\hbar^kQP,\] then
\[F*_{\hbar}G=Q_{\hbar}^{-1}(Q_{\hbar}F\circ Q_{\hbar}G),\;\;\forall F,G\in\op{Pol}(T^*M)\]
is such a $1$-parameter family.\\

As mentioned, the upshots of Sections \ref{derpoly} and
\ref{derdiffop} can be deduced from some invariant operators
obtained in \cite{NP} and \cite{BHMP}---although independent
proofs have also been found recently. The paper \cite{NP} fits
into the frame of a series of works, specially by P. Lecomte, P.
Mathonet, and E. Tousset \cite{LMT}, H. Gargoubi and V. Ovsienko
\cite{GO}, P. Cohen, Yu. Manin, and D. Zagier \cite{CMZ}, C. Duval
and V. Ovsienko \cite{DO}. A small dimensional hypothesis in
\cite{NP}, which was believed to be inherent in the used canonical
symbolization technique, was the starting point of \cite{BHMP}.
Here, the authors prove---using a conceptual method---existence
and uniqueness of a projectively equivariant symbol between the
space ${\cal D}^k_p(M)$ (${\cal D}^k_p$ for short) of $k$th order
differential operators transforming differential $p$-forms into
functions, and the space ${\cal S}^k_p(M)$ (${\cal S}^k_p$ for
short) of the corresponding symbols, the underlying manifold $M$
being endowed with a flat projective structure. This invariant
symbol map, which is explicitly known in terms of a ``divergence
operator'' and the ``Koszul differentials'', can be used as a
substitute for the affine or canonical symbol method and allows to
get rid of the dimensional assumption. Roughly speaking, the
search for ${\cal X}$-equivariant operators between ${\cal D}^k_p$
and ${\cal D}^{\E}_q$ reduces to investigations on
$\op{sl}$-equivariant operators on the classical level between
${\cal S}^k_p$ and ${\cal S}^{\E}_q$, and the subsequent
exploitation of the maximality of the projective algebra.

\section{Integrability of derivations}

Let us come back to the derivations of the algebras ${\cal
D}^1(M),$ ${\cal S}(M),$ and ${\cal D}(M)$. For these
non-integrable infinite-dimensional Lie algebras, there is no such
clear correspondence between derivations and $1$-parameter groups
of automorphisms as in the finite-dimensional setting. Our goal is
to find in each of these cases the most general form of a
$1$-parameter group of automorphisms. Moreover, computations
should unmask a derivation that can be viewed as the generator of
the chosen group of automorphisms. Finally, we wonder if it is
possible to characterize those derivations that induce
$1$-parameter groups of automorphisms.

First remark that any diffeomorphism $\zvf$ of $M$ canonically
induces an automorphism $\zvf_*$ of the considered algebra ${\cal
D}$. If ${\cal D}={\cal D}^1(M)$ or ${\cal D}={\cal D}(M)$, this
automorphism is defined by
\[(\phi_*D)f=D(f\circ\phi)\circ\phi^{-1},\;\;\forall D\in{\cal D},\forall f\in{\cal A}.\]
If ${\cal D}={\cal S}(M)= \op{Pol}(T^*M)$, we set
\[\phi_*P=P\circ(\phi^{\sharp})^{-1},\;\;\forall P\in{\cal D},\]
where $\phi^{\sharp}$ is the phase lift of $\zvf$. So 1-parameter
groups of diffeomorphisms are special 1-parameter groups of
automorphisms, known a priori, since they are just flows of
complete vector fields.

We have shown in \cite{GP} that any automorphism $\zF$ of ${\cal
D}^1(M)$ has the form
\begin{equation}\Phi(f+X)=\phi_*(X)+(K\,f+\zL\op{div}X+\zW(X))\circ\phi^{-1},\;\;\forall
f\in{\cal A},\forall X\in{\cal X},\label{autodiffop}\end{equation}
where $\phi\in\op{Diff}(M),$ $K\in\R\backslash\{0\},$ $\zL\in\R,$
and $\zW\in\zW^1(M)\cap\op{ker}\op{d}$ are uniquely determined by
the chosen automorphism. Let
$\Phi_t=\Phi_{\phi_t,K_t,\zL_t,\zW_t}$ be an arbitrary
$1$-parameter group of automorphisms. Smoothness with respect to
the differential structure of $M$ is assumed. In other words, we
suppose that the map \[\R\times M\ni
(t,x)\raa(\Phi_tD)(f)(x)\in\R\] is smooth for any $D\in{\cal
D}^1(M)$ and any $f\in{\cal A}$. When computing the l.h.s. of the
group condition
\begin{equation}\Phi_{\phi_t,K_t,\zL_t,\zW_t}\circ\Phi_{\phi_s,K_s,\zL_s,\zW_s}=\Phi_{\phi_{t+s},K_{t+s},\zL_{t+s},\zW_{t+s}},\label{groupcond}\end{equation}
we get terms that can easily be compared with the corresponding
terms of the r.h.s., except for one term,
\begin{equation}\zL_t\;(\op{div}\phi_{s*}X)\circ\phi_t^{-1},\label{badterm}\end{equation} which is not
of one of the four types in the r.h.s. of Equation
(\ref{autodiffop}). So this term has to be transformed.

Let us recall that the divergence is implemented by a fixed
nowhere vanishing $1$-density $\zr_0.$ It is quite obvious that
the divergence of the push-forward of a vector field $X$ coincides
with the divergence with respect to the pull-back of $\zr_0$. More
precisely,
\begin{equation}\op{div}_{\zr_0}\phi_*X=\lp\op{div}_{\phi^*\zr_0}X\rp\circ\phi^{-1},\label{pushpull}\end{equation}
where subscript $s$ has been omitted. It is clear that for any
diffeomorphism $\phi$ there is a unique positive smooth function
$J(\phi)$, such that \begin{equation}\phi^*\zr_0=\lp
J(\phi)\rp\zr_0.\label{Div}\end{equation} Furthermore, the reader
might have guessed that the essential local building block of
$J(\phi)(x)$ is $\mid\nolinebreak\op{det}\p_x\varphi\mid$, where
$\zf$ is the local form of $\phi$. Hence, the following property
of $J$:
\begin{equation}J(\phi\circ\psi)=\psi^*\lp J(\phi)\rp .\,
J(\psi),\;\;\forall\phi,\psi\in\op{Diff}(M).\label{precocycle}\end{equation}
Remember now that if $G$ is a group and $A$ is a left $G$-module,
a group 1-cocycle is a map $C:G\raa A$ such that
\begin{equation}C(g_1\cdot g_2)=g_1.\lp C(g_2)\rp+C(g_1),\;\;\forall g_1,g_2\in
G,\label{groupcocycle}\end{equation} where ``.'' is the action of
$G$ and ``$\cdot$'' the group multiplication. Note that this
cocycle-condition is similar to that of the Hochschild cohomology
of an associative algebra. The unique difference between Equation
(\ref{precocycle}) and Equation (\ref{groupcocycle}) is the
operation in the r.h.s. When applying the logarithm to both sides
of Equation (\ref{precocycle}), we finally get \[(\op{ln}\circ
J)(\phi\circ\psi)=\psi^*\lp (\op{ln}\circ J)(\phi)\rp+
(\op{ln}\circ J)(\psi).\] So \begin{equation}\op{ln}\circ
J\in{\cal Z}^1(\op{Diff}(M),\Ci(M))\label{Div2}\end{equation} is a
$1$-cocycle of the group of diffeomorphisms valued in the module
of smooth functions. Moreover, it can be proven that
\begin{equation}(\op{ln}\circ J)(\op{Exp}(tX))=\int_0^t\op{div}X\circ
\op{Exp}(sX)\op{ds},\label{Div3}\end{equation} for any complete
vector field $X$. Equations (\ref{Div}), (\ref{Div2}), and
(\ref{Div3}) show that \[\op{Div}=\op{ln}\circ J\] is the group
analogue of the divergence.

This group divergence allows to rewrite term (\ref{badterm}) in an
appropriate form. Starting from Equation (\ref{pushpull}), we
obtain
\[\begin{array}{ll}\lp\op{div}_{\phi^*\zr_0}X\rp\circ\phi^{-1}&=\lp\op{div}_{\lp J(\phi)\rp\zr_0}X\rp\circ\phi^{-1}\\&=\lp
\op{div}_{\zr_0}X+X\lp\lp\op{ln}\circ
J\rp(\zvf)\rp\rp\circ\zvf^{-1}\\&=\lp\op{div}\,X+\op{d}\lp
\op{Div}\zvf\rp(X)\rp\circ\zvf^{-1}.\end{array}\] The terms of
group condition (\ref{groupcond}) can now easily be compared. This
comparison leads to the equations
\[\begin{array}{c}\phi_t\circ\phi_s=\phi_{t+s},\phi_0=\op{id},\\
K_tK_s=K_{t+s},K_0=1,\\
\zL_t+K_t\zL_s,\zL_0=0,\\K_t\zW_s+\phi^*_s\zW_t+\zL_t\op{d}(\op{Div}\phi_s)=\zW_{t+s},\zW_0=0,\end{array}\]
where $\op{id}$ is the identity map. The first of these results
for instance means that the $1$-parameter family of
diffeomorphisms is actually a $1$-parameter group of
diffeomorphisms, i.e. the flow of a complete vector field $Y$:
$\phi_t=\op{Exp}(tY)$. Other equations are a little bit more
complicated, but can be solved, so that the explicit form of
$K_t,$ $\zL_t,$ and $\zW_t,$ i.e. of
$\Phi_{\phi_t,K_t,\zL_t,\zW_t}$ is known. Furthermore, the
solutions of the preceding equations involve, in addition to
vector field $Y$, two real numbers $\zk,\zl$ and a closed $1$-form
$\zw$. All these objects are uniquely determined by the chosen
$1$-parameter group of automorphisms and characterize, as
explained in Subsection \ref{derfristdiffop}, a derivation of
${\cal D}^1(M).$ We say that this derivation, which is special in
the sense that it is associated with a complete vector field,
induces the $1$-parameter group of automorphisms.

We are now ready to understand the following

\begin{theo} A derivation
\[C_{Y,\zk,\zl,\zw}(X+f)=[Y,X+f]+\kappa\,f+\lambda\,\op{div}X+\omega(X)\] of $\cD^1(M)$ induces
a one-parameter group $\zF_t$ of automorphisms of $\cD^1(M)$ if
and only if the vector field $Y$ is complete. In this case the
group is of the form
\beas&\Phi_t(X+f)=(\op{Exp}(tY))_*(X)+\left(e^{\zk
t}\,f+\,\zl\frac{e^{\zk t}-1}{\zk} \op{div}X\right)\circ
\op{Exp}(-tY)+\cr
&\left(\int_0^te^{\zk(t-s)}\left(\zl\int_0^sX(\op{div}Y\circ
\op{Exp}(uY))du+\lp(\op{Exp}(sY))^*\zw\rp(X)\right)ds\right)\circ
\op{Exp}(-tY).\label{1automd1} \eeas
\end{theo}

Similar upshots have been obtained for the algebras ${\cal S}(M)$
and ${\cal D}(M)$. They will not be described here.

\section{Differential operators on real line bundles}

In \cite{GP} we have taken an interest in characterizations of
manifold structures, especially by the Lie algebra of linear
differential operators acting on the functions of the chosen
manifold $M$. We now extend these results to the Lie algebra of
differential operators acting on tensor densities over $M$ of
arbitrary weight and even to differential operators acting on the
smooth sections of an arbitrary $\R$-line bundle $L$. Our
objectives are to examine if this Lie algebra structure recognizes
the base manifold $M$ and maybe even the bundle $L$ itself.

Let $\pi:L\raa M$ be a real vector bundle of rank 1 over a smooth,
Hausdorff, second countable, and connected manifold. We define the
algebra ${\cal D}(L)=\cup_{k\in\mathbb{N}}{\cal D}^k(L)$ of
differential operators on $L$ in the standard way. Note first that
the space $\op{Sec}(L)$ of smooth sections of $L$ is an $\cal
A$-module, so that any function $f\in{\cal A}$ induces an
endomorphism $m_f:\op{Sec}(L)\ni s\raa fs\in\op{Sec}(L)$ of the
space $\op{Sec}(L)$. Then set \[{\cal
D}^0(L)=\{D\in\op{End}(\op{Sec}(L)):[D,m_f]=0,\forall f\in{\cal
A}\},\]
\[{\cal D}^{k+1}(L)=\{D\in\op{End}(\op{Sec}(L)):[D,m_f]\in{\cal D}^k(L),
\forall f\in{\cal A}\}\;\;(k\in\mathbb{N}),\]
where $\op{End}(\op{Sec}(L))$ denotes the algebra of endomorphisms
of $\op{Sec}(L)$ and $[.,.]$-- the commutator bracket associated
with the composition multiplication.

\begin{prop} Any differential operator on $L$ is a local operator.\end{prop}

{\it Proof.} Indeed, if we denote ${\cal D}^{-1}(L)=\{ 0\}$, then
we can proceed now by induction and consider a $k$th-order ($k\ge
0$) differential operator $D$ on $L$, a section $s\in\op{Sec}(L)$
that vanishes in an open subset $U\subset M$, an arbitrary point
$x\in U$, and a function $\za\in{\cal A}$, which vanishes outside
$U$ and has constant value $1$ in some neighborhood of $x$ in $U$.
Since $[D,m_{\za}]\in{\cal D}^{k-1}(L)$ is a local operator,
$D(\za s)=\za D(s)$ in $U$, so $(D(s))(x)=D(0)=0$.
\rule{1.5mm}{2.5mm}

\begin{rem} In the following we write ${\cal D}(L)$ in the form ${\cal D}(L\raa M)$,
if we wish to put emphasis on the base manifold $M$, and in the
form ${\cal D}(M)$, if $L$ is the trivial bundle
$M\times\mathbb{R}$. This algebra ${\cal D}(M)$ is nothing but the
usual algebra of linear differential operators acting on the space
of smooth functions of $M$.
\end{rem}

\begin{prop}The space ${\cal D}(L)=\cup_{k\in\mathbb{N}}{\cal D}^k(L)$
is a quantum Poisson algebra in the sense of \cite{GP}.\end{prop}

{\it Proof.} In view of the above proposition it is sufficient to
check it locally. Let $U$ be an open subset of $M$ such that $L$
is trivial over $U$. Any local trivialisation of $L$ in $U$, i.e.
any nowhere vanishing section $\varsigma\in\op{Sec}(L_U)$ of $L$
over $U$, induces a canonical vector space isomorphism
\[\zi_{\varsigma}:\op{Sec}(L_U)\ni \zf\varsigma\raa \zf\in\Ci(U).\] This
isomorphism induces itself an isomorphism of quantum Poisson
algebras, \begin{equation}{\cal I}_{\varsigma}:{\cal D}(L_U)\ni
\zD\raa \zi_{\varsigma}\circ \zD\circ\zi_{\varsigma}^{-1}\in{\cal
D}(U).\label{identification}\end{equation} A straightforward
induction on the degree of differentiation allows to see that
${\cal I}_{\varsigma}$ respects the filtration.  \
\rule{1.5mm}{2.5mm}

\medskip\noindent
A gauge change entails a change in the identification of ${\cal
D}(L_U)$ with ${\cal D}(U).$ Indeed, if $\varsigma'$ is another
nowhere vanishing section of $L_U$, we have
$\varsigma'=\psi\varsigma$, $\psi\in\Ci(U)$, and, as easily
verified,
\begin{equation} {\cal I}_{\varsigma}(\zD)=m_{\psi}\circ{\cal
I}_{\varsigma'}(\zD)\circ m_{\psi^{-1}},\forall\zD\in{\cal
D}(L_U),\label{gaugetrans}\end{equation} so that ${\cal
I}_{\varsigma'}^{-1}\circ{\cal I}_{\varsigma}$ is a Lie algebra
automorphism of ${\cal D}(L_U)$. Hence the local isomorphisms
${\cal I}_{\varsigma}$ cannot be glued canonically. Note
nevertheless that if $\psi$ is a (non-zero) constant, the
identifications ${\cal I}_{\varsigma}$ and ${\cal I}_{\varsigma'}$
coincide.

Let us recall that quantum Poisson algebras (qPa) canonically
induce classical Poisson algebras (cPa) in the sense of \cite{GP}.
We denote by ${\cal S}(\cal D)$ the cPa implemented by a qPa
${\cal D}.$ We also know that ${\cal S}({\cal D}(M))$ coincides
with the algebra $\op{Sec}({\cal S}TM)$ of symmetric contravariant
tensor fields over $M$ and with the algebra $\op{Pol}(T^*M)$ of
smooth functions on the cotangent bundle of $M$ that are
polynomial along the fibers.

\begin{theo} The classical Poisson algebras ${\cal S}({\cal D}(L))$
and ${\cal S}({\cal D}(M))$ induced by ${\cal D}(L)$ and ${\cal
D}(M)$ are canonically isomorphic cPa.\end{theo}

{\it Proof.} In any qPa ${\cal D}$ we can define the $k$th-order
symbol $\zs_k(D)$ of any differential operator $D\in{\cal D}$, the
degree of which is $\le k$ (see \cite{GP}). In the fundamental
case ${\cal D}={\cal D}(M)$, this algebraically defined symbol
coincides with the usual geometric $k$th-order symbol. It is then
clear (see ($\ref{gaugetrans}$)) that \begin{equation}\zs_k({\cal
I}_{\varsigma}(\zD))=\zs_k({\cal I}_{\varsigma'}(\zD)),\forall
\zD\in{\cal D}^k(L_U).\label{global}\end{equation} It is obvious
that $\zs_k({\cal I}_{\varsigma}(\zD))\in{\cal S}_k({\cal D}(U))$,
$\zD\in{\cal D}^k(L_U)$ only depends on the $k$th-order symbol of
$\zD$, i.e. is actually defined on ${\cal S}_k({\cal D}(L_U))$.
Indeed, for any $P\in{\cal S}_k({\cal D}(L_U))$, if
$P=\zs_k(\zD)=\zs_k(\zD')$, we have $\zs_k({\cal
I}_{\varsigma}(\zD))-\zs_k({\cal I}_{\varsigma}(\zD'))=0$, since
$\zs_k(\zD-\zD')=0$. So \[\zF^U_{\varsigma}:{\cal S}_k({\cal
D}(L_U))\ni P\raa\zs_k({\cal
I}_{\varsigma}(\zs_k^{-1}(P)))\in{\cal S}_k({\cal D}(U))\] is a
cPa isomorphism. The morphism properties with respect to the
associative commutative and the Poisson-Lie multiplications are
direct consequences of the definitions of these operations $.$ and
$\{.,.\}$:
\[\zs(\zD).\zs(\zD')=\zs_{\op{deg}\zD+\op{deg}\zD'}(\zD\circ\zD'),\]
\[\{\zs(\zD),\zs(\zD')\}=\zs_{\op{deg}\zD+\op{deg}\zD'-1}([\zD,\zD']),\]
where $\zs(\zD)$ and $\op{deg}\zD$ are the principal symbol and
the degree of $\zD$ respectively (see \cite{GP}). In view of
(\ref{global}) the isomorphisms $\zF^{U}_{\varsigma}$ define a
global cPa isomorphism \[\zF:{\cal S}({\cal D}(L))\raa{\cal
S}({\cal D}(M)),\] such that \[(\zF
P)\m_U=\zF^{U}_{\varsigma}(P\m_U),\] for any $P\in{\cal S}_k({\cal
D}(L))$ and a trivialization $\varsigma$ of $L$ over a member $U$
of some appropriate open covering of $M$. Let us mention that if
$P=\zs_k(D)$, $D\in{\cal D}^k(L)$, the restriction
$P\m_U=(\zs_k(D))\m_U$ is nothing but the well-defined class
$\zs_k(D\m_U)$ of the restriction to $U$ of the local operator
$D$. \rule{1.5mm}{2.5mm}

\begin{theo} The quantum Poisson algebras ${\cal D}(L)$ and ${\cal D}(M)$
are isomorphic. \label{fundaresult}\end{theo}

{\it Proof.} Let $L_0=L\backslash \{0\}$ be the bundle $L$ with
removed $0$-section and let $\mid\! L_0\!\mid=L_0/\mathbb{Z}_2$ be
the quotient of $L_0$ with respect to the obvious action of the
multiplicative group $\mathbb{Z}_2=\{-1,1\}$. This quotient
$\mid\! L_0\!\mid$ is an affine bundle of dimension 1 (canonically
modelled on the vector bundle $M\times\mathbb{R}$), i.e. a smooth
bundle of 1-dimensional affine spaces, such that the passage from
one trivialization to another is given by an affine map. More
precisely, if $\varsigma$ is a never vanishing section of $L$ over
$U$, the trivialization of $\mid\!L_0\!\mid$ over $U$ is given by
\[\zf:U\times\mathbb{R}\ni (x,t)\raa
e^t\mid\!\varsigma(x)\!\mid\in\mid\!L_0\!\mid_U.\] Here
$\mid\!\varsigma(x)\!\mid$ denotes the class of $\varsigma(x)$ in
the quotient $\mid\!L_0\!\mid$ and the multiplication by $e^t$ is
the canonical multiplication of such a class by a non-zero real
number. A change of the non vanishing local section of $L$ is
characterized by a non vanishing smooth function
$\psi\in\Ci(U,\mathbb{R}^*)$, hence the corresponding change of
trivialization is given by the multiplication by
$\mid\!\psi\!\mid$, which is affine. The fibers of $\mid\!
L_0\!\mid$ are affine lines with a canonical free and transitive
action of $\mathbb{R}$ induced by the Liouville vector field of
$L$, $a_x.t=e^ta_x$ for any $t\in\mathbb{R}$ and any $a_x$ in the
fiber of $\mid\! L_0\!\mid$ over $x\in M.$ This action turns
$\mid\! L_0\!\mid$ into an $\mb{R}$-principal bundle. The
compatibility condition $\zf(x,s+t)=\zf(x,s).t$ $(x\in
U,s,t\in\mathbb{R})$ is obviously verified. Since the fibers are
contractible, $\mid\! L_0\!\mid$ has a global section
$\mid\!\zh\!\mid$. If $(U_{\za})_{\za\in\zL}$ is a covering of $M$
by open connected subsets over which $L$ is trivializable, section
$\mid\!\zh\!\mid$ can be viewed as a family
$(\{\zh_{\za},-\zh_{\za}\})_{\za\in\zL}$ of pairs of non vanishing
local sections of $L$, such that
$\{\zh_{\za},-\zh_{\za}\}=\{\zh_{\zb},-\zh_{\zb}\}$ on
$U_{\za}\cap U_{\zb}.$ This follows immediately from the above
depicted trivializations. When choosing for each $\za$ a
representative $\tilde{\zh}_{\za}\in\{\zh_{\za},-\zh_{\za}\}$, we
get a family of nowhere vanishing local sections of $L$, such that
on $U_{\za}\cap U_{\zb}$, we have $\tilde{\zh}_{\za}=\pm
\tilde{\zh}_{\zb}$ that reflexes the fact that line bundles over
$M$ are classified by $H^1(M;\Z_2)$. In view of
(\ref{identification}) and (\ref{gaugetrans}) we then get a global
quantum Poisson algebra isomorphism between ${\cal D}(L)$ and
${\cal D}(M)$.
\rule{1.5mm}{2.5mm}\\

The preceding result shows that the Lie algebra ${\cal D}(L\raa
M)$ characterizes the smooth structure of the base manifold $M$,
but does not recognize the topological complications in $L$.

\begin{cor}Let $\zp:L\raa M$ and $\zp':L'\raa M'$
be two real vector bundles of rank 1 over two smooth manifolds $M$
and $M'$ respectively. The Lie algebras ${\cal D}(L\raa M)$ and
${\cal D}(L'\raa M')$ are isomorphic if and only if the base
manifolds $M$ and $M'$ are
diffeomorphic.\label{charcresult}\end{cor}

{\it Proof.} Immediate consequence of Theorem \ref{fundaresult}
and of \cite[Theo. $6$]{GP}. \rule{1.5mm}{2.5mm}\\

Note that the essential fact is that ${\cal D}(L\raa M)$ and
${\cal D}(L'\raa M)$ are always isomorphic (even as qPa),
independently whether $L$ and $L'$ are isomorphic as vector
bundles or not. However these isomorphisms are not canonical ones,
depending on the choice of the section $\mid\!\zh\!\mid$ of
$\mid\! L_0\!\mid$. This observation, together with the
description of automorphisms and derivations of ${\cal D}(M)$ and
${\cal D}^1(M)$ \cite{GP,GP3}, gives automatically the obvious
description of automorphisms and derivations of ${\cal D}(L\raa
M)$ and ${\cal D}^1(L\raa M)$. The only difference is that
diffeomorphims $\phi$ of $M$ (resp. vector fields on $M$) do not
define automorphisms (resp. derivations) of ${\cal D}(L\raa M)$
and ${\cal D}^1(L\raa M)$ canonically, but in the way depending on
the choice of the section $\mid\!\zh\!\mid$ of $\mid\! L_0\!\mid$.
For example, the derivation $C_X^{\mid\!\zh\!\mid}$ associated
with a vector field $X$ on $M$ and a choice of $\mid\!\zh\!\mid$
is defined locally by
$$C_X^{\mid\!\zh\!\mid}(D)(f\zh_\alpha)
=[X(D_\alpha(f))-D_\alpha(X(f))]\zh_\alpha,$$ where
$D_\alpha(f)\zh_\alpha=D(f\zh_\alpha)$. This definition is
correct, since nothing changes when we choose $-\zh_\alpha$
instead of $\zh_\alpha$.

\noindent Janusz GRABOWSKI\\Polish Academy of Sciences\\Institute
of Mathematics\\\'Sniadeckich 8\\P.O. Box 21\\00-956 Warsaw,
Poland\\Email: jagrab@impan.gov.pl\\\\
\noindent Norbert PONCIN\\University of Luxembourg\\Mathematics
Laboratory\\avenue
de la Fa\"{\i}encerie, 162A\\
L-1511 Luxembourg City, Grand-Duchy of Luxembourg\\Email:
norbert.poncin@uni.lu

\end{document}